\documentclass[12pt,a4paper]{amsart}

\usepackage{latexsym} 
 
\usepackage[dvips]{graphics}
\usepackage{epsfig}

\usepackage{amssymb}
\usepackage{amsthm}
\usepackage{amsfonts}
\usepackage{amsmath}
\usepackage{amstext}
\usepackage{amscd}

\usepackage{epic}
\usepackage{eepic}
\usepackage{pstricks,pst-plot}

\numberwithin{equation}{section}

\theoremstyle{plain}
\newtheorem{thm}{Theorem}[section] 
\newtheorem{prop}[thm]{Proposition}
\newtheorem{cor}[thm]{Corollary}
\newtheorem{lem}[thm]{Lemma}

\newtheorem{theorem*}{Theorem}[]

\theoremstyle{definition}
\newtheorem{defn}[thm]{Definition}

\newtheorem{example}[thm]{Example}
\newtheorem{prob}{Problem}
\newtheorem{ques}{Question}

\theoremstyle{remark}
\newtheorem{rem}[thm]{Remark}

\newcommand{\N}{\mathbb{N}}
\newcommand{\R}{\mathbb{R}}

\def\accentclass@{7}
\def\makeacc@#1#2{\def#1{\mathaccent"\accentclass@#2 }}
\makeacc@\cir{017}

\title[the directional dimension of subanalytic sets]
{The  directional dimension  of subanalytic sets\\ 
is invariant
under bi-Lipschitz homeomorphisms}

\author{Satoshi Koike and Laurentiu Paunescu$^*$}
\address{Department of Mathematics, Hyogo University of Teacher Education,
Kato, Hyogo 673-1494, Japan}
\email{koike@sci.hyogo-u.ac.jp} 
\address{School of Mathematics, University of Sydney, Sydney, NSW, 2006,
Australia}
\email{laurent@maths.usyd.edu.au}

\date{}

\begin{document}

\maketitle


\begin{abstract}
Let $A \subset \R^n$ be a set-germ at $0 \in \R^n$
such that $0 \in \overline{A}$.
We say that $r \in S^{n-1}$ is a direction of $A$ at $0 \in \R^n$
if there is a sequence of points $\{ x_i \} \subset A \setminus \{ 0 \}$
tending to $0 \in \R^n$ such that ${x_i \over \| x_i \|} \to r$ 
as $i \to \infty$.
Let $D(A)$ denote the set of all directions of $A$ at $0 \in \R^n$.

Let $A, \ B \subset \R^n$ be subanalytic set-germs at $0 \in \R^n$
such that $0 \in \overline{A} \cap \overline{B}$.
We study the problem of whether the dimension of 
the common direction set, $\dim (D(A) \cap D(B))$
is preserved by  bi-Lipschitz homeomorphisms.
We show that although it is not true in general,  it is
preserved if the images of $A$ and $B$ 
are also subanalytic.
In particular if two subanalytic set-germs are  bi-Lipschitz
equivalent  their direction sets must have the same dimension.
\end{abstract}

\section{Introduction.}\label{introduction}

The first remarkable result on Lipschitz equisingularity problem was 
obtained by T. Mostowski. 
In \cite{mostowski1} he succeeded in solving  a conjecture of Sullivan,
showing that a complex analytic variety
admits a locally Lipschitz trivial stratification.
Following his work, A. Parusi\'nski proved the corresponding
results in several real categories (\cite{parusinski1,
parusinski2, parusinski3}).
Subsequently this area has become more attractive for real and 
complex singularity people.
Recently, J.P. Henry and A. Parusi\'nski (\cite{henryparusinski1,
henryparusinski2}) introduced some Lipschitz invariants for real and 
complex analytic function germs, and showed that Lipschitz moduli
appear even in a family of polynomial functions with isolated singularities. 
See the survey \cite{mostowski3} for more on Lipschitz 
equisingularity problems.

On the other hand, in late 70's, T.-C. Kuo introduced 
the notion of blow-analyticity as a desirable equivalence relation 
for real analytic function germs. He also  established some triviality 
theorems and showed local finiteness of different blow-analytic types in an
analytic family of functions with isolated singularities
(e.g. \cite{kuo3, kuo4, kuo5}).
Concerning blow-analyticity,
see the surveys \cite{fukuikoikekuo} and \cite{fukuipaunescu2}.

Let us recall the notion of blow-analyticity.
Let $f, \ g : (\R^n,0) \to (\R,0)$ be analytic function-germs.
We say that they are {\em blow-analytically equivalent}
if there are real modifications
$\mu : (M,\mu^{-1}(0)) \to (\R^n,0)$,
$\mu^{\prime} : (M^{\prime},\mu^{\prime -1}(0)) \to (\R^n,0)$
and an analytic isomorphism $\Phi : (M,\mu^{-1}(0)) \to
(M^{\prime},\mu^{\prime -1}(0))$ which induces a homeomorphism 
$\phi : (\R^n,0) \to (\R^n,0)$,  
$\mu^{\prime}\circ \Phi=\phi\circ \mu $, such that $f = g \circ \phi$.
A {\em blow-analytic homeomorphism} is such a $\phi$, a 
homeomorphism
induced by an analytic isomorphism via real modifications.

Every blow-analytic homeomorphism is an {\em arc-analytic} homeomorphism
in the sense of K. Kurdyka \cite{kurdyka},
therefore maps any analytic arc to an analytic arc. 
E. Bierstone and P. Milman analysed the relation between 
blow-analyticity and arc-analyticity in \cite{bierstonemilman}.
Taking those results into consideration, T.-C. Kuo conjectured that
a blow-analytic homeomorphism preserves the contact order 
of analytic arcs.
Nevertheless, this is not valid.
The first author observed that the zero-sets of  Brian\c{c}on-Speder's 
family (\cite{brianconspeder}) and also of Oka's family (\cite{oka})
are not ``blow-analytically and bi-Lipschitz'' trivial (in \cite{koike1},
see also \cite{paunescu}).
Later (in \cite{koike2}) he showed that 
they are not even bi-Lipschitz trivial
(while being blow-analytically trivial, see \cite{fukui},
\cite{fukuipaunescu1}).
In other words, the
blow-analytic equivalence for functions does not imply the
bi-Lipschitz equivalence for their zero-sets.
The proof in the case of Oka's family (see Example \ref{okafamily})
is based on the fact that
the cardinal number of the common direction set of their components
must be preserved by a bi-Lipschitz homeomorphism.
In this paper we extend the observation above
to the general case in the subanalytic category. 

\vspace{3mm}

\noindent {\bf Main Theorem.} {\em Let $A$, $B \subset \R^n$ be 
subanalytic set-germs at $0 \in \R^n$ such that 
$0 \in \overline{A} \cap \overline{B}$, and let 
$h : (\R^n,0) \to (\R^n,0)$ be a bi-Lipschitz homeomorphism.
Suppose that $h(A), \ h(B)$ are also subanalytic.
Then we have the equality of dimensions,}

\vspace{2mm}

\centerline{$\dim (D(h(A)) \cap D(h(B))) = \dim (D(A) \cap D(B)).$}

\vspace{4mm}

As a corollary of the theorem above, we have
another bi-Lipschitz invariant, namely the dimension of the direction set.

\vspace{3mm}

\begin{thm}\label{directionsetth} 
{\em Let $h : (\R^n,0) \to (\R^n,0)$ be 
a bi-Lipschitz homeomorphism. If $A, \ h(A)$ are subanalytic
set-germs at $0 \in \R^n$, then $\dim D(h(A)) = \dim D(A)$.} 
\end{thm}   
\vspace{3mm}

In \S 3 we describe our Main Problem and 
 give several examples showing the subtlety of our result. 
One example points out that the bi-Lipschitz assumption cannot be 
dropped, even if we 
deal with  polynomial homeomorphisms. 
Another two examples 
 demonstrate that we cannot drop the assumption of 
subanalyticity of the images from our main results. 
In \S 4 and \S 5 we define the notion of a sea-tangle
neighbourhood, describe some of its properties,
and introduce a sequence selection property (condition ($SSP$)).
After several reductions of our Main Problem in \S 6, 
we complete the proof in \S 7. 

At the end of this paper, we give an easy proof of the main theorem
for surfaces (see Appendix).

A special case of our result was obtained by Mostowski
in \cite{mostowski2}.


\vspace{0.2 truecm}  

{\bf Acknowledgements.} We would like to thank Tzee-Char Kuo
for useful discussions, in particular, for suggesting
the construction of a zigzag bi-Lipschitz homeomorphism
(Example \ref{zigzag}). 
This article was written up while the first author
was visiting Sydney.
He would like to thank the University of Sydney
for its support and hospitality. Finally we are happy to thank the referee for
helping improving the presentation of our paper.


\bigskip
\section{Directional dimension.}\label{kissing}
\medskip

We first recall the notion of subanalyticity introduced by
H. Hironaka (\cite{hironaka1}).
Let $M$ be a real analytic manifold.
A subset $A \subset M$ is said to be {\em subanalytic},
if for any $x \in \overline{A}$, there are an open neighbourhood
$U$ of $x$ in $M$ and a finite numbers of proper real analytic maps 
of real analytic spaces $f_{ij} : Y_{ij} \to U$, $j =1,2$, such that
$$
A \cap U = \bigcup_i \ (Im(f_{i1}) - Im(f_{i2})).
$$
There are several equivalent definitions for subanalyticity
(\cite{hironaka1, hironaka2}).
We note that the curve selection lemma, called Hironaka's selection 
lemma, holds in the subanalytic category.

We next give the definition of the direction set.

\begin{defn}\label{directionset}
Let $A$ be a set-germ at $0 \in \R^n$ such that
$0 \in \overline{A}$.
We define the {\em direction set} $D(A)$ of $A$ at $0 \in \R^n$ by
$$
D(A) := \{a \in S^{n-1} \ | \
\exists  \{ x_i \} \subset A \setminus \{ 0 \} ,
\ x_i \to 0 \in \R^n  \ \text{s.t.} \
{x_i \over \| x_i \| } \to a, \ i \to \infty \}.
$$
Here $S^{n-1}$ denotes the unit sphere centred at $0 \in \R^n$.
\end{defn}

\noindent Thanks to Hironaka's selection lemma,
we can express the direction set $D(A)$ for a subanalytic set-germ
$A$ at $0 \in \R^n$ as follows:
$$
D(A) := \biggl\{a \in S^{n-1} \ | \
\begin{matrix}
\exists  \lambda : [0,\epsilon ) \to \R^n, \ C^{\omega}, \ 
\lambda (0) = 0,\  \lambda ((0,\epsilon )) \setminus \{ 0 \} 
\subset A \\
\ \text{s.t.} \ 
\lim_{t \to 0} {\lambda^{\prime}(t) \over \| \lambda^{\prime}(t) \|} = a 
\end{matrix}\biggr\}.
$$

Concerning this direction set, we have

\begin{prop}\label{direction}
If $A$ is a subanalytic set-germ at $0 \in \R^n$ such that
$0 \in \overline{A}$, then $D(A)$ is a closed subanalytic subset
of $S^{n-1}$.   
\end{prop} 

\begin{proof}
Let $\pi : \mathcal{M}_n \to \R^n$ be a blowing-up at $0 \in \R^n$
such that $\pi^{-1}(0) = \R P^{n-1}$.
Let $\beta : S^{n-1} \to \R P^{n-1}$ be the canonical projection,
and we write $\hat{P} := \beta (P)$ for $P \in S^{n-1}$.

Let $\epsilon > 0$ be a fixed sufficiently small positive number.
For $Q \in \R P^{n-1}$, we denote by $U_{\epsilon}(Q)$
the $\epsilon$-neighbourhood of $Q$ in $\mathcal{M}_n$.
Then $U_{\epsilon}(Q) - \pi^{-1}(0) = U_{\epsilon}^+(Q) \cup
U_{\epsilon}^-(Q)$, where $U_{\epsilon}^+(Q), \ U_{\epsilon}^-(Q)$
are disjoint open half balls.

We denote by $T$ the strict transform of $\overline{A}$ by $\pi$.
Let $P$ be an arbitrary point of $S^{n-1}$.
Then there exists a neighbourhood $U$ of $P$ in $S^{n-1}$ such
that $D(A) \cap U$ can be identified with $\pi^{-1}(0) \cap 
\overline{T \cap U_{\epsilon}^+(\hat{P})} \cap U_{\epsilon}(\hat{P})$
or $\pi^{-1}(0) \cap 
\overline{T \cap U_{\epsilon}^-(\hat{P})} \cap U_{\epsilon}(\hat{P})$,
which is a closed subanalytic set in $U_{\epsilon}(\hat{P})$.
Thus $D(A)$ is a closed subanalytic subset of $S^{n-1}$.
\end{proof}

Let $A$, $B$ be subanalytic set-germs at $0 \in \R^n$
such that $0 \in \overline{A} \cap \overline{B}$.
By the proposition above, $D(A) \cap D(B)$
is a closed subanalytic subset of $S^{n-1}$.
Therefore the dimension of $D(A) \cap D(B)$ is naturally defined
(by convention $\dim \emptyset = -1$).

\begin{defn}\label{kissingdimension}
For subanalytic set-germs $A$, $B$ at $0 \in \R^n$
such that $0 \in \overline{A} \cap \overline{B}$,
we call $\dim (D(A) \cap D(B))$ the {\em directional dimension}
of $A$ and $B$ at $0 \in \R^n$.
\end{defn}

\begin{rem}\label{remark1}
Let $A \subset \R^n$ be a subanalytic set-germ at $0 \in \R^n$
such that $0 \in \overline{A}$, and let $h : (\R^n,0) \to (\R^n,0)$
be a bi-Lipschitz homeomorphism.
Since a  subanalytic subset of $\R^n$ admits a locally finite 
stratification by connected analytic submanifolds of $\R^n$,
$h(A)$ admits a finite stratification by connected Lipschitz
submanifolds of $\R^n$ and $\dim h(A) = \dim A$.
\end{rem}

Let us apply our Main Theorem to Oka's family (\cite{oka}).

\begin{example}\label{okafamily}
Let 
$f_t : (\R^3,0) \to (\R,0)$, $t \in \R$,
be a family of polynomial functions
with isolated singularities defined by
$$
f_t(x,y,z) = x^8 + y^{16} + z^{16} + t x^5 z^2 + x^3 y z^3.
$$
We recall some observations in \cite{koike1}. 
Put
$$
f(x,y,z) := f_0(x,y,z) = x^8 + y^{16} + z^{16} + x^3 y z^3.
$$
The set $f^{-1}(0) - \{ 0 \}$ has empty intersection
with each coordinate plane.
Let us consider

\vspace{3mm}

\qquad $A_1 := \{ x > 0, \ y > 0, \ z < 0 \}, \ \ 
A_2 := \{ x > 0, \ y < 0, \ z > 0 \}$,

\qquad $A_3 := \{ x < 0, \ y > 0, \ z > 0 \}, \ \ 
A_4 := \{ x < 0, \ y < 0, \ z < 0 \}$

\vspace{3mm}

\noindent and $S_i := f^{-1}(0) \cap A_i$, $1 \le i \le 4$.
Then $f^{-1}(0) = S_1 \cup S_2 \cup S_3 \cup S_4 \cup \{ 0 \}$
and each $\overline{S_i} = S_i \cup \{ 0 \}$ is 
homeomorphic to $S^2$.
As seen in \cite{koike1}, $\dim (D(S_i) \cap D(S_j)) = 0$,
$i \ne j$.

We further introduce

\vspace{3mm}

\qquad $A_5 := \{ x < 0, \ y < 0, \ z > 0 \}, \ \ 
A_6 := \{ x < 0, \ y > 0, \ z < 0 \}$.

\vspace{3mm}

\noindent The zero-set $f_t^{-1}(0)$ is expanding into the octants 
$A_5$ and $A_6$ as $t$ varies from $0$ to $1$.
In \cite{koike1}, we have made the following observation for $f_1^{-1}(0)$.
Put
$$
g(x,y,z) := f_1(x,y,z) = x^8 + y^{16} + z^{16} + x^5 z^2 + x^3 y z^3.
$$
The set $g^{-1}(0) - \{ 0 \}$ has empty intersection with
both $(x,y)$-plane and $(y,z)$-plane.
We put

\vspace{3mm}

\qquad $B_1 := \{ x > 0, \ y > 0, \ z < 0 \}, \ \ 
B_2 := \{ x > 0, \ y < 0, \ z > 0 \}$,

\qquad $B_3 := \{ x < 0, \ z > 0 \}, \ \ 
B_4 := \{ x < 0, \ z < 0 \}$

\vspace{3mm}

\noindent and $P_i := g^{-1}(0) \cap B_i$, $1 \le i \le 4$.
Then $g^{-1}(0) = P_1 \cup P_2 \cup P_3 \cup P_4 \cup \{ 0 \}$
and each $\overline{P_i} = P_i \cup \{ 0 \}$ is 
homeomorphic to $S^2$.
We have seen $\dim (D(P_3) \cap D(P_4)) = 1$.
Thus it follows from our Main Theorem that $(\R^3,f_0^{-1}(0))$ 
is not bi-Lipschitz homeomorphic to $(\R^3,f_1^{-1}(0))$. 
In fact, 
 the same argument  shows that the zero sets of $f_0$ and $f_t, t\neq 0$ are 
not bi-Lipschitz homeomorphic.
\end{example}


\bigskip
\section{Main problem and examples of bi-Lipschitz
homeomorphisms.}\label{sea-tanglenbd}
\medskip

Here we pose the following natural question: 

\vspace{3mm}

\noindent {\bf Main Problem}
\begin{prob}\label{mainprob} 
 {\em Let $A$, $B \subset \R^n$ be 
subanalytic set-germs at $0 \in \R^n$ such that 
$0 \in \overline{A} \cap \overline{B}$, and let 
$h : (\R^n,0) \to (\R^n,0)$ be a bi-Lipschiptz homeomorphism.
Suppose that $h(A), \ h(B)$ are also subanalytic.
Then is it true that}
\end{prob}
\vspace{2mm}

\centerline{$\dim (D(h(A)) \cap D(h(B))) = \dim (D(A) \cap D(B))$ ?}

\vspace{5mm}
The next example points out that the bi-Lipschitz assumption cannot be 
dropped, even if we 
deal with  polynomial homeomorphisms. 
\begin{example}\label{anahomeo}
Let $h :  (\R^3,0) \to (\R^3,0)$ be the polynomial homeomorphism  defined by
 $h(x,y,z) = (x,y,z^3)$.
 The variety $V=\{x^2+y^2-z^6=0\}$ is mapped onto the
variety $W=\{x^2+y^2-z^2=0\}$. Clearly they have different directional 
dimensions.
\end{example}

We now offer two examples of bi-Lipschitz homeomorphisms which demonstrate
that we cannot drop the assumption that the images are also subanalytic.

\begin{example}\label{spiral} (Quick spiral). 
Let $h = (h_1,h_2) : (\R^2,0) \to (\R^2,0)$ be a map defined by

\vspace{3mm}

\qquad $h_1(x,y) = x \cos (\log (x^2+y^2)) + y \sin (\log (x^2+y^2))$,

\vspace{2mm} 

\qquad $h_2(x,y) = - x \sin (\log (x^2+y^2)) + y \cos (\log (x^2+y^2))$,\

\vspace{3mm}

\noindent in other words, $h(r,\theta ) = (r,\theta - \log r)$
in the polar coordinates. A half-line with the initial point
at the origin is mapped by $h$ to a spiral below: 

\vspace{3mm}

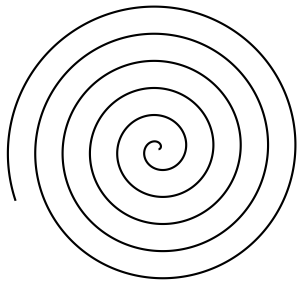
\begin{figure}[ht]
\begin{pspicture}*(-2,-2)(2,2)
  \parametricplot[plotpoints=2000]{0}{2000}
  {t cos t mul .001 mul t sin t mul .001 mul}
\end{pspicture}
\caption{A spiral}
\end{figure}

\vspace{3mm}

\noindent Then it is easy to see that ${\partial h_1 \over \partial x}$,
${\partial h_1 \over \partial y}$, ${\partial h_2 \over \partial x}$,
${\partial h_2 \over \partial y}$ are bounded in a punctured
neighbourhood of $0 \in \R^2$.
Therefore $h$ is Lipschitz near $0 \in \R^2$.
Similarly, we can see that $h^{-1}$ is also Lipschitz.
Thus $h$ is a bi-Lipschitz homeomorphism.

Let $A$, $B$ be two different segments with an end point at $0 \in \R^2$.
Then their images have $D(h(A)) = D(h(B)) = S^1$, which implies 
$\dim (D(h(A)) \cap D(h(B))) = 1$. 
But, it is clear that $D(A) \cap D(B) = \emptyset$, which
implies $\dim (D(A) \cap D(B)) = -1$.
\end{example}

\begin{example}\label{zigzag} (Zigzag bi-Lipschitz homeomorphism).
Let $f : (\R,0) \to (\R,0)$ be a zigzag function whose graph is
drawn below, namely $f$ is zigzag if $x \in (0,1)$,
and $f \equiv 0$ if $x \notin (0,1)$. For instance we can take 
$a_n:=(\frac{\sqrt 3 -1}{\sqrt 3 +1})^n, n$ a non-negative integer, and define 
 $f(x) = \sqrt 3 (x-a_n)$ if $x\in [a_n,a_{n-1}\frac{\sqrt 3} {1+\sqrt 3}]$, 
and  $f(x) = \sqrt 3 (-x+a_{n-1})$ if 
$x\in [a_{n-1}\frac{\sqrt 3} {1+\sqrt 3},a_{n-1} ]$.

\vspace{3mm}

\epsfxsize=4cm
\epsfysize=3cm
\begin{figure}[htb]
$$\epsfbox{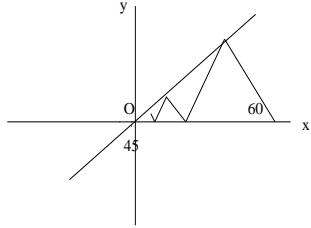}$$
\caption{zigzag function}
\label{fig:zigzag}
\end{figure}

\vspace{3mm}

\noindent It is easy to see that $f$ is a Lipschitz function. 
Define a map $h = (h_1,h_2) : (\R^2,0) \to (\R^2,0)$ by
$$
h_1(x,y) = x, \ \ h_2(x,y) = y + f(x).
$$
Then $h$ is a bi-Lipschitz homeomorphism.

Let $A$, $B$ be two different segments with an end point at $0 \in \R^2$ 
and such that $A$ is on the positive $x$-axis, 
and $B$ is very close to $A$, namely the angle at the origin  
between them is very small.
Then we can see that $\dim (D(h(A)) \cap D(h(B))) = 1$, but
$\dim (D(A) \cap D(B)) = -1$.
\end{example}


\bigskip
\section{Sea-tangle neighbourhood and properties.}\label{sea-tanglepro}
\medskip

In this section we define the notion of a sea-tangle
neighbourhood for a subset of $\R^n$.

\begin{defn}\label{seatanglenbd}
Let $A \subset \R^n$ such that $0 \in \overline{A}$,
and let $d, \ C > 0$.
The {\em sea-tangle neighbourhood $ST_d(A;C)$ of $A$,
of degree $d$ and  width} $C$, is defined by:
$$
ST_d(A;C) := \{ x \in \R^n \ | \ dist (x,A) \le C |x|^d \}.
$$
\end{defn}
This definition originated from the classical notion of 
{\em horn-neighbourhood} (e.g. T.C. Kuo \cite{kuo1, kuo2}).
In fact, if $A$ is an analytic arc $ST_d(A;C)$ is horn-like;
if $A$ is a tangling Lipschitz arc it looks like a sea-tangle.

Let $\mathcal{S}$ be the set of set-germs $A \subset \R^n$
at $0 \in \R^n$ such that $0 \in \overline{A}$. 
We next introduce an equivalence relation in $\mathcal{S}$.

\begin{defn}\label{STequivalence}
Let $A$, $B \in \mathcal{S}$.
We say that $A$ and $B$ are $ST$-{\em equivalent},
if there are $d_1$, $d_2 > 1$, $C_1$, $C_2 > 0$ such that
$B \subset ST_{d_1}(A;C_1)$ and $A \subset ST_{d_2}(B;C_2)$
as germs at $0 \in \R^n$. 
We write $A \sim_{ST} B$.
\end{defn}

\begin{rem}\label{remark2}
It is easy to see that the $ST$-equivalence $\sim_{ST}$ is an equivalence 
relation in $\mathcal{S}$.
\end{rem}

Let $\phi : (\R^n,0) \to (\R^n,0)$ be a bi-Lipschitz homeomorphism,
namely there are positive numbers $K_1, \ K_2 > 0$
with $K_1 \le K_2$ such that
$$
K_1 |x_1 - x_2| \le |\phi (x_1) - \phi (x_2)| 
\le K_2 |x_1 - x_2|
$$
in a small neighbourhood of $0 \in \R^n$.
Conversely, we have
$$
\frac{1}{K_2} |y_1 - y_2| \le |\phi^{-1}(y_1) - \phi^{-1}(y_2)| 
\le \frac{1}{K_1} |y_1 - y_2|
$$
in a small neighbourhood of $0 \in \R^n$.
In \cite{koike2}, we have shown that a kind of Sandwich Lemma holds
for the sea-tangle neighbourhoods of a Lipschitz arc and of its image 
by a bi-Lipschitz homeomorphism.
Using a similar argument, we can show the following:

\begin{lem}\label{sandwichlemma}(Sandwich Lemma).
Let $A \subset \R^n$ such that $0 \in \overline{A}$.
Then, for $K > 0$, 
$$
ST_d(\phi(A);\frac{KK_1}{K_2^d}) \subset \phi(ST_d(A;K)) \subset
ST_d(\phi(A);\frac{KK_2}{K_1^d})
$$
in a small neighbourhood of $0 \in \R^n$.
\end{lem}

By this Sandwich Lemma, we can easily see the following proposition:

\begin{prop}\label{remark3}
$ST$-equivalence is preserved by a bi-Lipschitz homeomorphism.
\end{prop}

We introduce some notations.
For a subset $A \subset S^{n-1}$, we denote by $L(A)$
a half-cone of $A$ with the origin $0 \in \R^n$ as the vertex:
$$
L(A) := \{ t a \in \R^n\ | \ a \in A, \ t \ge 0 \}.
$$
We make some notational conventions.
In the case $A = \{ a\}$, we simply write $L(a) := L(\{ a \} )$.
For a set-germ $A$ at $0 \in \R^n$ such that $0 \in \overline{A}$, 
we put $LD(A) := L(D(A))$, the {\em real tangent cone}
at $0 \in \R^n$. 

\begin{example}\label{curve}
Let $\pi : \mathcal{M}_2 \to \R^2$ be a blowing-up at $(0,0) \in \R^2$,
and let $a = (0,1) \in S^1$.
We denote by $\hat{L}(a)$ the strict transform of $L(a)$ 
in $\mathcal{M}_2$ by $\pi$.
In a suitable coordinate neighbourhood,
$\pi : \R^2_{(X,Y)} \to \R^2$ can be expressed as $\pi (X,Y) = (XY,Y)$.
Here $(0,0) \in \R^2_{(X,Y)}$ is the intersection of $\hat{L}(a)$
and the exceptional divisor $E = \pi^{-1}(0,0)$.

Let
$B := \{ (X,Y)\in \R^2_{(X,Y)} \ | \ Y = e^{- {1 \over |X|^2}},
\ X \ge 0 \}$.
Then the curve $B$ is not contained in
$\{ (X,Y) \in \R^2_{(X,Y)} \ | \ 
|Y| \ge C^{\prime} |X|^{d^{\prime}} \}$
as germs at $(0,0) \in \R^2_{(X,Y)}$,
for any $d^{\prime} > 0$, $C^{\prime} > 0$.

Let $A := \pi (B)$. 
Then we can see that $\lim_{m \to \infty} {a_m \over \| a_m \|} = a$
for any sequence of points $\{ a_m \}$ on $A$ tending to
$(0,0) \in \R^2$, which implies $LD(A) = L(a)$.
Moreover $A$ is not contained
in any sea-tangle neighbourhood $ST_d(LD(A);C)$
as germs at $(0,0) \in \R^2$, for $d > 1$, $C > 0$. 
\end{example}

On the other hand, in the subanalytic case we have the following:

\begin{prop}\label{keyproperty}
Let $A$ be a subanalytic set-germ at $0 \in \R^n$ such that 
$0 \in \overline{A}$. 
Then there is $d_1 > 1$ such that $A \subset ST_d(LD(A);C)$ as set-germs 
at $0 \in \R^n$ for any $d$ with $1 < d  < d_1$ and $C > 0$.
\end{prop}

\begin{proof} 
Since the order of $d(\gamma(t), LD(A))$ is greater than the order of 
$\gamma(t)$ on each analytic arc at $0$ in $A$, 
the function $g(x) = {d(x,LD(A)) \over \| x \|}$ extends 
at the origin as $g(0) = 0$ (use Hironaka's selection lemma).
The Lojasiewicz inequality (\cite{lojasiewicz}, \cite{bochankrisler})
for $g(x)$ and $\| x \|$ gives that
$g(x) \le \| x \|^{\epsilon}$, for some $\epsilon > 0$,
in a small neighbourhood of $0 \in \R^n$.
Setting $d_1 = 1 + \epsilon > 1$,
the statement holds for any $d$ with $1 < d  < d_1$, $C > 0$.
\end{proof}

We next describe the key lemma for analytic arcs;
it takes an important role in the proof of our Appendix.
We denote by $\mathcal{A}(\R^n,0)$ the set of germs of analytic 
maps $\lambda : [0,\epsilon ) \to \R^n$ with
$\lambda (0) = 0, \ \lambda (s) \ne 0, \ s > 0$.
For any $\lambda \in \mathcal{A}(\R^n,0)$, there exists a unique
$a \in S^{n-1}$ such that $\lambda$ is tangent to $L(a)$
at $0 \in \R^n$.
Then we write $T(\lambda ) := L(a)$.

\begin{lem}\label{keyspecial} (Key Lemma for analytic arcs).
Let $h : (\R^n,0) \to (\R^n,0)$ be a bi-Lipschitz homeomorphism.
Suppose that there are $\gamma_1, \ \gamma_2 \in \mathcal{A}(\R^n,0)$
such that $T(\gamma_1 ) = T(\gamma_2 )$.
Then for any sequence of points $\{ a_m \} \subset h(\gamma_1 )$
tending to $0 \in \R^n$
with $\lim_{m \to \infty} {a_m \over \| a_m \|} = a \in S^{n-1}$,
there is a sequence of points $\{ b_m \} \subset h(\gamma_2 )$
tending to $0 \in \R^n$
such that $\lim_{m \to \infty} {b_m \over \| b_m \|} = a$
$($i.e. $D(h(\gamma_1)) = D(h(\gamma_2)))$.
\end{lem}

\begin{proof}
Since $T(\gamma_1 ) = T(\gamma_2 )$, there are $d > 1, C_1 > 0$
such that $\gamma_1 \subset ST_d(\gamma_2;C_1)$ as germs at $0 \in \R^n$.
By Lemma \ref{sandwichlemma}, there is $C_2 > 0$ such that
$h(\gamma_1) \subset ST_d(h(\gamma_2);C_2)$ as germs at $0 \in \R^n$.
Therefore, for any sequence of points $\{ a_m \} \subset h(\gamma_1 )$
tending to $0 \in \R^n$
with $\lim_{m \to \infty} {a_m \over \| a_m \|} = a$,
$B_{C_2 \| a_m \|^d}(a_m) \cap h(\gamma_2) \ne \emptyset$ for any $m$.
Here $B_r(P)$ denotes a ball centred at $P \in \R^n$ of radius $r > 0$.
For each $m$, take $b_m$ from the above intersection.
Let $\{ b_k \}$ be an arbitrary subsequence of $\{ b_m \}$ such that
$\lim_{k \to \infty} {b_k \over \| b_k \|} = b \in S^{n-1}$.
Suppose that $b \ne a$. 
Then there is $C_3 > 0$ such that
$$
ST_1(L(a);2C_3) \cap ST_1(L(b);2C_3) = \{ 0 \}.
$$
If $k$ is sufficiently large, we can assume that
$a_k \in ST_1(L(a);C_3)$, $\ b_k \in ST_1(L(b);C_3)$.
But $b_k \in B_{C_2 \| a_k \|^d}(a_k)$ implies
$b_k \in ST_1(L(a);2C_3)$ for sufficiently large $k$,
since $d > 1$.
This is a contradiction.
Thus $b = a$. 
\end{proof}

Now we discuss some sea-tangle properties in a more general setup.
Throughout this section,
let $A$, $B \subset \R^n$ be set-germs at $0 \in \R^n$
such that $0 \in \overline{A} \cap \overline{B}$,
namely $A$, $B \in \mathcal{S}$,
and let $h : (\R^n,0) \to (\R^n,0)$ be a bi-Lipschitz homeomorphism.
Then we can rewrite Lemma \ref{keyspecial} in the following form:

\begin{lem}\label{keygeneral} (Key Lemma for general sets).
Suppose that there are $d > 1$, $C > 0$ such that 
$A \subset ST_d(B;C)$ as germs at $0 \in \R^n$.
Then we have $D(h(A)) \subset D(h(B))$.
In addition, we have $D(ST_d(h(A));C^{\prime})) \subset D(h(B))$
for any $C^{\prime} > 0$.
\end{lem}

We have some corollaries of this lemma.

\begin{cor}\label{cor1}
$D(ST_d(h(A);C)) = D(h(A))$ for any $d > 1$, $C > 0$.
\end{cor}

\begin{cor}\label{cor3}
$D(ST_d(A;C)) = D(A)$ for any $d > 1$, $C > 0$.
\end{cor}

\begin{cor}\label{cor2}
Suppose that there are $d > 1$, $C > 0$ such that 
$A \subset ST_d(B;C)$ as germs at $0 \in \R^n$.
Then we have $D(A) \subset D(B)$.
In particular, if $A$ and $B$ are $ST$-equivalent,
then we have $D(A) = D(B)$.
\end{cor}
  
In the subanalytic case we give more sea-tangle properties.

\begin{prop}\label{oppositekey}
Suppose that $A$ is subanalytic. 
Then, for $d_1 > 1$, $C_1 > 0$, there is $1 < d_2 < d_1$
such that $ST_{d_1}(LD(A);C_1) \subset ST_d(A;C)$
as germs at $0 \in \R^n$, for any $d$ with $1 < d < d_2$ and $C > 0$.
\end{prop}

\begin{proof}
Since $d(\gamma(t), LD(A))\leq C_1 \| \gamma(t)\|^{d_1}$
on each analytic arc at $0$ contained in $ST_{d_1}(LD(A);C_1)$, 
we have that the order of $d(\gamma(t), A)$ is greater than the order 
of $\gamma(t)$. 
Using the same arguments as in Proposition \ref{keyproperty}
we conclude that for all $x\in ST_{d_1}(LD(A);C_1)$
we have $d(x,A)\leq |x|^{d_2}$ 
for some $d_2$ with $1 < d_2 < d_1$.
Therefore the statement holds for any $d$ with $1 < d < d_2$ and $C > 0$.
\end{proof}

The assumption of subanalyticity is essential in 
Proposition \ref{oppositekey}. 
For instance, see Example \ref{zigzag}.

By Propositions \ref{keyproperty}, \ref{oppositekey}, we have

\begin{thm}\label{eqtheorem} 
If $A$ is subanalytic, then $A$ is $ST$-equivalent to $LD(A)$.
\end{thm}

As a corollary of Proposition \ref{oppositekey}, we have

\begin{cor}\label{cor4}
Suppose that $h(A)$, $h(B)$ are subanalytic.
If $D(h(A)) \subset D(h(B))$, then  there are $d > 1$, $C > 0$
such that $A \subset ST_d(B;C)$ as germs at $0 \in \R^n$.
\end{cor}

\begin{proof}
By Proposition \ref{keyproperty} and the assumption, 
there are $d_1 > 1$, $C_1 > 0$ such that
$$
h(A) \subset ST_{d_1}(LD(h(A));C_1) \subset ST_{d_1}(LD(h(B));C_1)
$$
as germs at $0 \in \R^n$.
By Proposition \ref{oppositekey},
there are $1 < d < d_1$, $C_2 > 0$ such that
$$
ST_{d_1}(LD(h(B));C_1) \subset ST_d(h(B);C_2)
$$
as germs at $0 \in \R^n$.
Thus we have 
$$
h(A) \subset ST_d(h(B);C_2)
$$
as germs at $0 \in \R^n$.
Then it follows that
$$
A = h^{-1}(h(A)) \subset h^{-1}(ST_d(h(B);C_2))
$$
as germs at $0 \in \R^n$.
By Lemma \ref{sandwichlemma}, there is $C > 0$ such that
$$
A \subset h^{-1}(ST_d(h(B);C_2)) \subset ST_d(B;C)
$$
as germs at $0 \in \R^n$.
\end{proof}

Using the results above we can characterise
the conditions in the Key Lemma as follows:

\begin{thm}\label{iff}
Suppose that $h(A)$, $h(B)$ are subanalytic.
Then the following conditions are equivalent.

(1) $D(h(A)) \subset D(h(B))$.

(2) There are $d > 1$, $C > 0$ such that $A \subset ST_d(B;C)$
as germs at $0 \in \R^n$.
\end{thm}


\bigskip
\section{Sequence selection property.}\label{Condition ($SSP$)}
\medskip

In this section we introduce a sequence selection property,
and discuss some consequences for  the sets satisfying it .

\begin{defn}\label{conditionT}
Let $A \subset \R^n$ be a set-germ at $0 \in \R^n$
such that $0 \in \overline{A}$.
We say that $A$ satisfies {\em condition} $(SSP)$,
if for any sequence of points $\{ a_m \}$ of $\R^n$
tending to $0 \in \R^n$ such that 
$\lim_{m \to \infty} {a_m \over \| a_m \| } \in D(A)$,
there is a sequence of points $\{ b_m \} \subset A$ such that
$$
\| a_m - b_m \| \ll \| a_m \|, \ \| b_m \| . 
$$ 
\end{defn}

\begin{example}\label{exampleT}
(1) Let $A := \{ b_m \} \subset \R$ be a sequence of points
defined by
$$
b_{m+1} = (1 - 2 \epsilon )b_m, \ \ 0 < \epsilon < {1 \over 2},
$$
where $b_1 > 0$. 
Let $a_m := (1 - \epsilon )b_m$, \ $m \in \N$.
Then $a_m = {b_m + b_{m+1} \over 2}$. 
Therefore we have
$$
D(\{ b_m \} ) = D(\{ a_m \} ) = \{ 1 \} \ \ \text{and}
\ \ |a_m - b_m| = |a_m - b_{m+1}| =  \epsilon |b_m|.
$$
Thus $A$ does not satisfy condition $(SSP)$.

Let $B := \{ b_m \} \subset \R$ be a sequence of points
defined by $b_m = {1 \over m}$. 
Then $B$ satisfies condition $(SSP)$.

(2) Let $T$ be an angle with vertex at $O\in \R^2$.
We choose sequences of points $\{ P_m \}$ and $\{ Q_m \}$ 
on the edges of $T$ such that $\overline{OP_m} = {1 \over m}$ and 
$\overline{OQ_m}$ has its abscisa ${1 \over 2}({1 \over m} + {1 \over m+1})$
(see the figure below).
Let $C_1$ be a zigzag curve connecting $P_m$'s and $Q_m$'s.

\epsfxsize=4cm
\epsfysize=3cm
\begin{figure}[htb]
$$\epsfbox{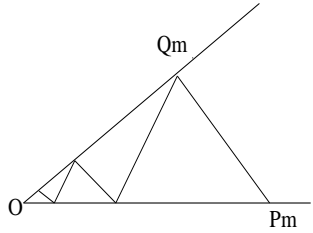}$$
\caption{zigzag curve}
\label{fig:zigzag}
\end{figure}

\vspace{3mm}

\noindent Then $C_1$ satisfies condition $(SSP)$.
Since the length of $C_1$ is infinite,
$C_1$ is not an image of any subanalytic curve by any bi-Lipschitz 
homeomorphism. 

If instead we choose $\{ P_m \}$ and $\{ Q_m \}$
such that $\overline{OP_m} = {1 \over m^2}$ and 
$\overline{OQ_m} = {1 \over 2}({1 \over m^2} + {1 \over (m+1)^2})$
and $C_2$ is a zigzag curve connecting $P_m$'s and $Q_m$'s,
$C_2$ satisfies condition $(SSP)$ and the length of $C_2$ is finite.

(3) The curve $A$ defined in Example \ref{curve}
satisfies condition $(SSP)$.
\end{example}

We make some remarks.

\begin{rem}\label{remark4}
Let $A \subset \R^n$ be a set-germ at $0 \in \R^n$
such that $0 \in \overline{A}$.

(1) The cone $LD(A)$ satisfies condition $(SSP)$.

(2) If $A$ is subanalytic, then it satisfies condition $(SSP)$.
\end{rem}

\begin{rem}\label{remark5}
Condition $(SSP)$ is $C^1$ invariant
but not bi-Lipschitz invariant (see Example \ref{zigzag} ).
\end{rem} 

\begin{rem}\label{remark6}
We would like to emphasise the fact that
$A \sim_{ST} LD(A)$ is specific to the subanalytic category.
If $A$ satisfies merely condition $(SSP)$,
this does not always guarantee that $A \sim_{ST} LD(A)$ 
(see Examples \ref{curve} and \ref{exampleT} (3)).
\end{rem}

As in the previous section,
let $A$, $B \subset \R^n$ be set-germs at $0 \in \R^n$
such that $0 \in \overline{A} \cap \overline{B}$,
and let $h : (\R^n,0) \to (\R^n,0)$ be a bi-Lipschitz homeomorphism.
Here we show an important lemma, necessary for the proof of our Main Theorem.

\begin{lem}\label{directionrel}
$D(h(A)) \subset D(h(LD(A))).$
Moreover, if $A$ satisfies condition $(SSP)$, then the equality holds.
\end{lem}

\begin{proof}
For any $\alpha \in D(h(A))$, there is a sequence of points
$\{ a_m \} \subset A$ tending to $0 \in \R^n$ such that
$\lim_{m \to \infty} {h(a_m) \over \| h(a_m) \| } = \alpha$.
Then there is a subsequence $\{ a_k \}$ of $\{ a_m \}$ such that
$$
\lim_{k \to \infty} {a_k \over \| a_k \| } \in D(A) = D(LD(A)).
$$
Since $LD(A)$ satisfies condition $(SSP)$, there is a sequence of points 
$\{ b_k \} \subset LD(A)$ such that 
$\| a_k - b_k \| \ll \| a_k\| , \ \| b_k\|$.
It follows that 
$\|h(a_k)- h(b_k)\| \ll \| h(a_k)\|, \ \| h(b_k)\|$.
Thus
$$
\alpha = \lim_{k \to \infty} {h(a_k) \over \| h(a_k) \| }
=\lim_{k \to \infty} {h(b_k) \over \| h(b_k) \| } \in D(h(LD(A))),
$$
that is $D(h(A)) \subset D(h(LD(A)))$.

By replacing $A$ by $LD(A)$, we can similarly show the equality part.
\end{proof}

As a corollary of this lemma we have

\begin{cor}\label{direction1} $D(A) \subset D(h^{-1}(LD(h(A)))$.
\end{cor}

Using a similar argument as in Lemma \ref{directionrel},
we can show the following:

\begin{prop}\label{direction2}
Suppose that $B$ satisfies condition $(SSP)$.
If $D(A) \subset D(B)$, then $D(h(A)) \subset D(h(B))$.
\end{prop}

As a corollary of this proposition we have the following theorem:

\begin{thm}\label{control}
Suppose that $B$, $h(B)$ satisfy condition $(SSP)$.
Then $D(A) \subset D(B)$ if and only if $D(h(A)) \subset D(h(B))$.
\end{thm}

It is natural to ask the following question: 

\begin{ques}\label{question1}
Suppose that $A$, $B$ are subanalytic.
Then $D(A) \subset D(B)$ if and only if $D(h(A)) \subset D(h(B))$?
\end{ques}

The answer to this question is ``no''.
The ``if'' part does not always hold.
See Example \ref{spiral}.


\bigskip
\section{Reductions of Main Problem.}\label{reduction}
\medskip

Let $h : (\R^n,0) \to (\R^n,0)$ be a bi-Lipschitz homeomorphism,
and let $A_1$, $A_2$, $B_1$, $B_2 \subset \R^n$ be 
subanalytic set-germs at $0 \in \R^n$ such that
$0\in \overline{A_1} \cap \overline{A_2}$,
$0\in \overline{B_1} \cap \overline{B_2}$
and $h(A_i) = B_i$, $i = 1, 2$.

\vspace{3mm}

Let $A \subset \R^n$ be a set-germ at $0 \in \R^n$
such that $0 \in \overline{A}$.
Here we consider the following problem:

\begin{prob}\label{problem2}\
Suppose that $A$, $h(A)$ are subanalytic.
Then  is it true that
$$
\dim D(A) \ge \dim D(h(A)) \ ?
$$
\end{prob}

\begin{rem}\label{reamark7}
If the answer to Problem \ref{problem2} is affirmative,
then we have
$$
\dim D(A) = \dim D(h(A)).
$$
\end{rem}

\noindent Concerning this problem we have the following statement:

\vspace{3mm}

\noindent {\bf Statement.} We can reduce our Main Problem 
\ref{mainprob}
to Problem \ref{problem2}.

\begin{proof}
Indeed suppose that the answer to Problem \ref{problem2}
is affirmative.
Using Corollary \ref{cor3}, we can easily show the following equality:
$$ 
D(A_1) \cap D(A_2) =D(ST_{d_1}(LD(A_1);C_1) \cap ST_{d_2}(LD(A_2);C_2))
$$
for $d_i > 1$, $C_i > 0$, $i = 1, \ 2$. Therefore we have
$$
\dim (D(A_1) \cap D(A_2)) = 
\dim D(ST_{d_1}(LD(A_1);C_1) \cap ST_{d_2}(LD(A_2);C_2)).
$$
Since $A_1$, $A_2$ are subanalytic, by Theorem \ref{eqtheorem},
this also equals to
$$
\dim D(ST_{d_1}(A_1;C_1^{\prime}) \cap ST_{d_2}(A_2;C_2^{\prime}))
$$
for some $C_1^{\prime}$, $C_2^{\prime} > 0$.
Then it follows from Problem \ref{problem2} and Lemma 
\ref{sandwichlemma} that
$$
\dim D(ST_{d_1}(A_1;C_1^{\prime}) \cap ST_{d_2}(A_2;C_2^{\prime})) =
\dim D(ST_{d_1}(B_1;K_1) \cap ST_{d_2}(B_2;K_2))
$$
for some $K_1$, $K_2 > 0$.
Since the latter dimension equals to $\dim (D(B_1) \cap D(B_2))$,
we have
$$
\dim (D(A_1) \cap D(A_2)) = \dim (D(B_1) \cap D(B_2)).
$$
\end{proof}

\begin{rem}
Suppose that $A$, $h(A)$ are subanalytic. Then
$$
\dim LD(A)=\dim h(LD(A))\,\, \text{and} \dim LD(h(LD(A))) = \dim LD(h(A)) 
$$ 
so Problem \ref{problem2} is equivalent to showing that 
$$
\dim h(LD(A)) \ge \dim LD(h(LD(A))) 
$$

\end{rem}

 The remark above will give us 
the possibility to replace $A$ by its cone 
$LD(A)$ whenever convenient. Although $h(LD(A))$ is not subanalytic in 
general, it is more than just merely an image of a subanalytic set
by a bi-Lipschitz homeomorphism, it satisfies condition $(SSP)$.
In order to see this fact, we mention a lemma without proof.

\begin{lem}\label{seqlemma}
Let $A \subset \R^n$ be
 a set-germ at $0 \in \R^n$ such that 
$0 \in \overline{A}$, and let $d > 1$, $C > 0$.
For any sequence of points $\{ b_m \} \subset ST_d(A;C)$
tending to $0 \in \R^n$, there is a sequence of points
$\{ a_m \} \subset A$ such that $\| a_m - b_m \| \ll \| a_m \|^{d_1}$
for any $d_1$ with $1 \le d_1 < d$.
\end{lem}

\begin{prop}\label{h(LD(A))}
The set $h(LD(A))$ satisfies condition $(SSP)$.
\end{prop}

\begin{proof}
Let $\{ a_m \}$ be an arbitrary sequence of points of $\R^n$
tending to $0 \in \R^n$ such that
$$
\lim_{m \to \infty} {a_m \over \| a_m \| } \in D(h(LD(A)))
= D(h(A)).
$$
Since $h(A)$ is subanalytic, there is a sequence of points
$\{ b_m \} \subset h(A)$ such that
$$
\| a_m - b_m \| \ll \| a_m \|, \ \|b_m \| .
$$
This implies $\| b_m \| \le 2 \| a_m \|$ for sufficiently large $m$.
Since $A$ is also subanalytic, it follows from Proposition \ref{keyproperty}
that there are $d > 1$, $C > 0$ such that
$A \subset ST_d(LD(A);C)$ as germs at $0 \in \R^n$.
By Lemma \ref{sandwichlemma}, there is $C_1 > 0$ such that
$h(A) \subset ST_d(h(LD(A));C_1)$ as germs at $0 \in \R^n$.
It follows that $\{ b_m \} \subset ST_d(h(LD(A));C_1)$.
Then, by Lemma \ref{seqlemma}, there is a sequence of points
$\{ c_m \} \subset h(LD(A))$ such that
$\| c_m - b_m \| \ll \| c_m \|^{d_1}$ for any $d_1$ with $1 \le d_1 < d$.
This implies $\| c_m - b_m \| \ll \| c_m \|, \ \| b_m \|$
and $\| b_m \| \le 2\| c_m \|$ for sufficiently large $m$.
Therefore we have
$$
\| a_m - c_m \| \le \| a_m - b_m \| + \| b_m - c_m \| \ll 
\| b_m \|, \ \| a_m \| , \ \| c_m \| .
$$ 
Thus $h(LD(A))$ satisfies condition $(SSP)$.
\end{proof}


\bigskip
\section{Proof of main results.}\label{proof}
\medskip

We first make an observation on the volume of sea-tangle
neighbourhoods.

\begin{lem}\label{volumelemma} Let $\alpha$, $\beta$ 
be linear subspaces of $\R^n$.
Suppose that $\dim \alpha < \dim \beta$.
Then, for $d > 1$, $C_1$, $C_2 > 0$,
$$
\lim_{\epsilon \to 0} { Vol(ST_d(\alpha;C_1) \cap B_{\epsilon}(0))
\over Vol(ST_d(\beta;C_2) \cap B_{\epsilon}(0)) } = 0.
$$
\end{lem}

\begin{proof}
Put 
$$
\Gamma = \Gamma_{\alpha ,\beta} := \{ \tilde{\alpha}, \
\text{vector subspace of} \ \R^n \ | \ \tilde{\alpha} \subset \beta, \
\dim \tilde{\alpha} = \dim \alpha \}.
$$
Fix $C > 0$ and take $\epsilon > 0$.
For each $\tilde{\alpha} \in \Gamma$,
define $A_{\tilde{\alpha}} := ST_d(\tilde{\alpha};C) \cap B_{\epsilon}(0)$.
Let $\mu_{\epsilon}$ be the greatest number of pairwise
disjoint $A_{\tilde{\alpha}}$, $\tilde{\alpha} \in \Gamma$
such that $A_{\tilde{\alpha}} \subset ST_d(\beta ;C) \cap B_{\epsilon}(0)$.
Note that this number is necessarily finite.

Since $\mu_{\epsilon}$ tends to $\infty$ as $\epsilon \to 0$,
it follows that
$$
\lim_{\epsilon \to 0} { Vol(ST_d(\alpha;C) \cap B_{\epsilon}(0))
\over Vol(ST_d(\beta;C) \cap B_{\epsilon}(0)) } = 0.
$$
The fact that 
$$
Vol(ST_d(\alpha ;C_1) \cap B_{\epsilon}(0)) \le 
K \ Vol(ST_d(\alpha ;C) \cap B_{\epsilon}(0))
$$
where $K := ({C_1 \over C})^{n - \dim \alpha}$, implies
our observation.
\end{proof}

This lemma suggests that the same volume property
holds for the cones of subananlytic set-germs,
since a subanalytic set of $\R^n$ admits a locally finite stratification
by analytic submanifolds of $\R^n$ which are analytically 
equivalent to Euclidean spaces.

Let $f, \ g : [0,\delta ) \to \R$, $\delta > 0$, be non-negative
functions. 
If there are $K > 0$, $0 < \delta_1 \le \delta$ such that
$$
f(\epsilon ) \le K g(\epsilon ) \ \ \text{for} \ \ 0 \le \epsilon \le \delta_1,
$$
then we write $f \precsim g$ (or $g \succsim f$).
If $f \precsim g$ and $f \succsim g$, we write $f \thickapprox g$.

\begin{prop}\label{conevolume}
Let $\alpha$, $\beta \subset \R^n$ be subanalytic cones at $0 \in\R^n$ .
Suppose that $\dim \alpha < \dim \beta$.
Then, for $d > 1$, $C_1$, $C_2 > 0$,
$$
\lim_{\epsilon \to 0} { Vol(ST_d(\alpha;C_1) \cap B_{\epsilon}(0))
\over Vol(ST_d(\beta;C_2) \cap B_{\epsilon}(0)) } = 0.
$$
\end{prop}

\begin{proof}
Let $\gamma$ be a subanalytic cone at $0 \in \R^n$ of dimension $r$,
and let $M$ be an $r$-dimensional linear subspace of $\R^n$.
Then the proposition follows easily from Lemma \ref{volumelemma}
and the fact that
$$
Vol(ST_d(\gamma ;C) \cap B_{\epsilon}(0)) \thickapprox 
Vol(ST_d(M;C) \cap B_{\epsilon}(0))
$$
for $d > 1$, $C > 0$.
To see this fact, one may assume that $\gamma$ is equidimensional.
In this case we have
$$
ST_d(\gamma ;C) \subset \bigcup ST_d(T_x ;C),
$$
where the union is finite and $T_x$, $x \in \overline{\gamma} \cap S^{n-1}$, 
is an $r$-dimensional linear subspace of $\R^n$ through $x$.
This implies 
$$
Vol(ST_d(\gamma ;C) \cap B_{\epsilon}(0)) \precsim
Vol(ST_d(M;C) \cap B_{\epsilon}(0)).
$$

On the other hand, for  $x \in \gamma \cap S^{n-1}$,  
$\gamma$ is locally bi-Lipschitz equivalent to
the tangent space $T_x$ of $\gamma$ at $x$.
For $C, \delta > 0$, there is $K > 0$ such that
$$
Vol(ST_d(T_x \cap L(\tilde{B}_x(\delta ));C) \cap B_{\epsilon}(0)) 
\ge K \ Vol(ST_d(T_x;C) \cap B_{\epsilon}(0))
$$
for any small $\epsilon > 0$, where $\tilde{B}_x(\delta )$ is 
a $\delta$-neighbourhood of $x$ in $S^{n-1}$.
Thus we can claim the opposite inequality $\succsim$ as well.   
\end{proof}

In general, we have the following relation on dimensions
for subanalytic set-germs: 

\begin{lem}\label{subanalyticdimension}
Let $A \subset \R^n$ be a subanalytic set-germ at $0 \in \R^n$
such that $0 \in \overline{A}$.
Then we have $\dim LD(A) \le \dim A$.
\end{lem}

\begin{proof}
Let $f : \overline{A} - \{ 0 \} \to S^{n-1}$ be the mapping defined by
$f(a) = {a \over \| a \|}$, and let
$\pi : \overline{\text{Graph} f} \to \R^n$ be the canonical projection.
Then $D(A) = D(\overline{A}) = \pi^{-1}(0)$.
Therefore we have
$$
\dim D(A) = \dim \pi^{-1}(0) < \dim \text{Graph} f = \dim \overline{A}
= \dim A.
$$
Thus it follows that $\dim LD(A) = \dim D(A) + 1 \le \dim A$.
\end{proof}

In addition, we have the following volume property on $ST$-equivalence:

\begin{prop}\label{STvolume}
Let $A, B \subset \R^n$ be set-germs at $0 \in \R^n$
such that $0 \in \overline{A} \cap \overline{B}$.
Suppose that $A$ and $B$ are $ST$-equivalent.
Then for $C_1, \ C_2 > 0$, there is $d_1 > 1$ such that
$$
Vol(ST_d(A;C_1) \cap B_{\epsilon}(0)) \thickapprox
Vol(ST_d(B;C_2) \cap B_{\epsilon}(0))
$$
for any $d$ with $1 < d \le d_1$.
\end{prop}

\begin{proof}
Since $A$ and $B$ are $ST$-equivalent, there are $d_3, \ d_4 > 1$
and $C_3, \ C_4 > 0$ such that $A \subset ST_{d_3}(B;C_3)$
and $B \subset ST_{d_4}(A;C_4)$ as germs at $0 \in \R^n$.
Let $d_1 = \min (d_3,d_4) > 1$.
Then for any $d$ with $1 < d \le d_1$, we have
$$
ST_d(A;C_1) \subset ST_d(ST_{d_3}(B;C_3);C_1) \subset ST_d(B;C_5)
$$
as germs at $0 \in \R^n$, where $C_5 = C_1 + C_3 > 0$.
Note that there is $K > 0$ such that
$$
Vol(ST_d(B;C_5) \cap B_{\epsilon}(0)) \le
K \ Vol(ST_d(B;C_2) \cap B_{\epsilon}(0))
$$
for any small $\epsilon > 0$. It follows that
$$
Vol(ST_d(A;C_1) \cap B_{\epsilon}(0)) \precsim 
Vol(ST_d(B;C_2) \cap B_{\epsilon}(0)).
$$
The opposite inequality $\succsim$ follows similarly.
\end{proof}

The following corollary is an obvious consequence of
Theorem \ref{eqtheorem}, Lemma \ref{subanalyticdimension} and
Propositions \ref{conevolume}, \ref{STvolume}.

\begin{cor}\label{volumeratio}
Let $\alpha \subset\R^n$ be a subanalytic set-germ at $0 \in\R^n$ 
such that $0 \in \overline{\alpha}$,
and let $\beta \subset \R^n$ be a subanalytic cone at $0 \in \R^n$.
Suppose that $\dim \alpha < \dim \beta$.
Then, for $C_1$, $C_2 > 0$, there is $d_1 > 1$ such that
$$
\lim_{\epsilon \to 0} { Vol(ST_d(\alpha;C_1) \cap B_{\epsilon}(0))
\over Vol(ST_d(\beta;C_2) \cap B_{\epsilon}(0)) } = 0.
$$
for any $d$ with $1 < d \le d_1$.
\end{cor}

\begin{rem}\label{remark8}
We cannot take $\beta$ merely a subanalytic set-germ
in the corollary above.
Let $\alpha \subset \R^3$ be the positive $z$-axis,
and let $\beta := \{ (x,y,z) \in \R^3 \ | \ z^3 = x^2 + y^2 \}$.
Then $\dim \alpha = \dim LD(\beta ) = 1$ and $\dim \beta = 2$.
For $d > 1$ sufficiently close to $1$ and $C > 0$,
$$
\lim_{\epsilon \to 0} { Vol(ST_d(\alpha;C) \cap B_{\epsilon}(0))
\over Vol(ST_d(\beta;C) \cap B_{\epsilon}(0)) } = 1.
$$
\end{rem}

Using Corollary \ref{volumeratio}, we can show the following lemma:

\begin{lem}\label{maintool}
Let $h : (\R^n,0) \to (\R^n,0)$ be a bi-Lipschitz homeomorphism,
let $E \subset \R^n$ be a subanalytic set-germ at $0 \in \R^n$
such that $0 \in \overline{E}$, and let $F := h(E)$.
Suppose that $F$ and $LD(F)$ are $ST$-equivalent and
$LD(F)$ is subanalytic.
Then we have $\dim LD(F) \le \dim E$.
\end{lem}

\begin{proof}
Assume that $\dim LD(F) > \dim F \ (= \dim E)$.
Since $F$ and $LD(F)$ are $ST$-equivalent,
it follows from Proposition \ref{STvolume} that
there are $d_1 > 1$ and $C_1, \ C_2 > 0$ such that
$$
Vol(ST_d(F;C_1) \cap B_{\epsilon}(0)) \thickapprox
Vol(ST_d(LD(F);C_2) \cap B_{\epsilon}(0))
$$
for any $d$ with $1 < d \le d_1$.

On the other hand, $h$ is a bi-Lipschitz homeomorphism.
Therefore we have the following volume relation:
$$
Vol(ST_d(F;C_1) \cap B_{\epsilon}(0)) \thickapprox  
Vol(ST_d(E;C_3) \cap B_{\epsilon}(0))
$$
for $C_3 > 0$. It follows that
$$
1 \thickapprox {Vol(ST_d(F;C_1) \cap B_{\epsilon}(0)) \over
Vol(ST_d(LD(F);C_2) \cap B_{\epsilon}(0))} \thickapprox  
{Vol(ST_d(E;C_3) \cap B_{\epsilon}(0)) \over
Vol(ST_d(LD(F);C_2) \cap B_{\epsilon}(0))}
$$
for $d$ with $1 < d \le d_1$.
By Corollary \ref{volumeratio}, the right ratio tends to $0$ as 
$\epsilon \to 0$, if $d > 1$ is sufficiently close to $1$.
This is a contradiction.
Thus we have $\dim LD(F) \le \dim F$.
\end{proof}

Now we show our Main Theorem.
By the reduction of Main Problem in the previous section,
it suffices to show that the answer to Problem \ref{problem2} is affirmative.
Let us recall the hypotheses of Problem \ref{problem2},
namely $h : (\R^n,0) \to (\R^n,0)$ is a bi-Lipschitz homeomorphism
and $A$, $h(A) \subset \R^n$ are subanalytic set-germs at $0 \in \R^n$
such that $0 \in \overline{A}$.

We apply Lemma \ref{maintool}  to
 $E := LD(A)$ and $F := h(LD(A))$, so we need to check all the assumptions 
of \ref{maintool}.

Because   $ h(A)$ is assumed subanalytic, so it is 
$LD(h(A))=LD(h(LD(A)))=LD(F)$.

Since $A$ is subanalytic, $LD(A)$ is $ST$-equivalent
 to $A$ (see Theorem \ref{eqtheorem}).
 Then, by Proposition \ref{remark3},
 $F=h(LD(A))$ is $ST$-equivalent to $h(A)$.
 In addition, it follows from the subanalyticity of $h(A)$ that
 $h(A)$ is $ST$-equivalent to $LD(h(A)) = LD((h(LD(A)))=LD(F)$.
Since $ST$-equivalence is an equivalence relation (Remark \ref{remark2}),
$F=h(LD(A))$ is $ST$-equivalent to $LD(F)=LD(h(LD(A)))$.

\vspace{2mm}

\noindent Therefore it follows from Lemma \ref{maintool} 
that $\dim LD(h(A))=\dim LD(h(LD(A)))  \le \dim LD(A)$,
which proves that the answer to Problem \ref{problem2} is affirmative,
and as a result, our Main Problem has an affirmative answer as well.
This  concludes the proof of our Main Theorem.

Obviously our Main Theorem can be generalized to arbitrary finite families of 
subanalytic sets.

Since we have shown  the affirmative answer to Problem \ref{problem2},
we have  proved Theorem \ref{directionsetth} as well, which also 
follows as a corollary of our  
Main Theorem.

\medskip

\begin{rem} The authors are preparing a  note with Ta L\^e Loi
on  directional properties in  $o$-minimal structures.
In that note we are also  discussing  whether the main result of this paper 
holds
in a $o$-minimal structure, replacing the assumptions of
subanalytic sets with those of definable sets.
The main result holds in a $o$-minimal structure over the real
field. However the natural correpoding result does not always hold
in a $o$-minimal structure over a general real closed field.
In fact the direction set can be infinite-dimensional.
In addtion, we used the finite covering property of compactness
(bounded closed sets) in our volume arguments, but compactness
does not mean the finite covering property over a general real closed field.
\end{rem}
\bigskip
\centerline{\bf Appendix.}
\medskip

In this appendix we give a quick proof of our Main Theorem
for subanalytic surfaces.
Let $f : (\R^n,0) \to (\R^p,0)$ be a subanalytic map-germ
such that $f^{-1}(0) - \{ 0 \} \ne \emptyset$ as germs at $0 \in \R^n$.
Then, for two connected components $A_1$, $A_2$ of
$f^{-1}(0) - \{ 0 \}$ (if they exist), 
$\overline{A_1} \cap \overline{A_2} = \{ 0 \}$.
Therefore we consider our Main Problem in the following setup:

Let $h : (\R^n,0) \to (\R^n,0)$ be a bi-Lipschitz homeomorphism,
and let $A_1$, $A_2$, $B_1$, $B_2 \subset \R^n$ be 
subanalytic set-germs at $0 \in \R^n$ such that
$\overline{A_1} \cap \overline{A_2} = \{ 0 \}$,
$\overline{B_1} \cap \overline{B_2} = \{ 0 \}$
and $h(A_i) = B_i$, $i = 1, 2$.

Under this setup we have the following claim on the directional dimension:

\vspace{3mm}

\noindent {\bf Claim 1.} $\dim (D(A_1) \cap D(A_2)) \le n - 2$ \ \
($\dim (D(B_1) \cap D(B_2)) \le n - 2$).

\begin{proof}
Since $D(A_1)$, $D(A_2) \subset S^{n-1}$,
we have $\dim D(A_1)$, $\dim D(A_2) \le n - 1$.
Suppose that

\vspace{3mm}

\noindent (A.1)  \ \ \ \ \ \ \ \ \ \ \ \ \ \ \ \ \ \ \ \ \ \ \
 $\dim (D(A_1) \cap D(A_2)) = n - 1$.

\vspace{3mm}

\noindent Then $\dim D(A_1) = \dim D(A_2) = n - 1$.
It follows from Lemma \ref{subanalyticdimension} that

\vspace{3mm}

\noindent (A.2)  \ \ \ \ \ \ \ \ \ \ \ \ \ \ \ \ \ \ \ \ \ \ \
$\dim A_1 = \dim A_2 = n$.

\vspace{3mm}

\noindent Then, by (A.1) and (A.2), 
$(A_1 - \{ 0 \}) \cap (A_2 - \{ 0 \}) \ne \emptyset$
as germs at $0 \in \R^n$, which contradicts our assumption.
Therefore we have  $\dim (D(A_1) \cap D(A_2)) \le n - 2$.
\end{proof}

By Lemma \ref{keyspecial}, we have 

\vspace{3mm}

\noindent {\bf Claim 2.} {\em If} $\dim (D(A_1) \cap D(A_2)) = -1$,
{\em then} $\dim (D(B_1) \cap D(B_2)) = -1$.

\vspace{3mm}

As seen in Proposition \ref{direction}, $D(A_1)$, $D(A_2)$ and
$D(A_1) \cap D(A_2)$ are closed subanalytic subsets of $S^{n-1}$.
Therefore they are compact.
In particular, if their dimension is $0$, they are finite points sets.

Concerning the directional dimension, we have another claim.

\vspace{3mm}

\noindent {\bf Claim 3.} {\em If} $\dim (D(A_1) \cap D(A_2)) = 0$,
{\em then} $\dim (D(B_1) \cap D(B_2)) = 0$.

\begin{proof}
If $\dim (D(B_1) \cap D(B_2)) = -1$, then, by Claim 2,
$\dim (D(A_1) \cap D(A_2)) = -1$.
Therefore $\dim (D(B_1) \cap D(B_2)) \ge 0$.

Suppose that

\vspace{3mm}

\noindent (A.3)  \ \ \ \ \ \ \ \ \ \ \ \ \ \ \ \ \ \ \ \ \ \ \
$\dim (D(B_1) \cap D(B_2)) \ge 1$.

\vspace{3mm}

\noindent Since $\dim (D(A_1) \cap D(A_2)) = 0$,
$D(A_1) \cap D(A_2)$ is a finite points set.
Let $D(A_1) \cap D(A_2) := \{ P_1, \cdots , P_a \}$ 
where $1 \le a < \infty$.
By (A.3), we can pick up $a+1$ points $Q_1, \cdots , Q_{a+1}$
from a connected subanalytic subset of $D(B_1) \cap D(B_2)$
of dimension $\ge 1$.
Corresponding to each $Q_j$, $1 \le j \le a + 1$,
there are analytic arcs $\alpha_j \subset B_1 \cup \{ 0 \}$, 
$\beta_j \subset B_2 \cup \{ 0 \}$ such that
$T(\alpha_j) = T(\beta_j) = L(Q_j)$.
Then it follows from Lemma \ref{keyspecial} that 
for any sequence of points $\{ a_m \} \subset h^{-1}(\alpha_j)$
such that $\lim_{m \to \infty} {a_m \over \| a_m \|}$ exists,
there exists a sequence of points $\{ b_m \} \subset h^{-1}(\beta_j)$
such that $\lim_{m \to \infty} {b_m \over \| b_m \|} =
\lim_{m \to \infty} {a_m \over \| a_m \|}$, $ 1 \le j \le a + 1$.

Here we make a remark on the limit point set.

\vspace{3mm}

\noindent {\em Remark} A. For each $j$, if   
$\lim_{m \to \infty} {a_m \over \| a_m \|}$ and
$\lim_{m \to \infty} {a_m^{\prime} \over \| a_m^{\prime} \|}$
exist for $\{ a_m \}$, $\{ a_m^{\prime} \} \subset h^{-1}(\alpha_j)$,
then their limit points coincide.
After this, we denote by $R_j$ the unique limit point.

\begin{proof}
Suppose that
$$
\lim_{m \to \infty} {a_m \over \| a_m \|} = a \ne a^{\prime}
= \lim_{m \to \infty} {a_m^{\prime} \over \| a_m^{\prime} \|}.
$$
Let $S_{\epsilon}(a) := L(\partial (B_{\epsilon}(a) \cap S^{n-1}))$.
Then there are $\epsilon_1$, $\epsilon_2 > 0$ with
$0 < \epsilon_1 < \epsilon_2 < \| a - a^{\prime} \|$
such that for any $\epsilon$ with $\epsilon_1 \le \epsilon \le \epsilon_2$,
$S_{\epsilon}(a) \cap h^{-1}(\alpha_j)$ contains infinitely
many points $\{ C_k^{\epsilon} \}$.
Therefore, for any $\epsilon$ with $\epsilon_1 \le \epsilon \le \epsilon_2$,
there is a subsequence $\{ C_t^{\epsilon} \}$ of $\{ C_k^{\epsilon} \}$
such that $\lim_{t \to \infty} {C_t^{\epsilon} \over \| C_t^{\epsilon} \| }
= C^{\epsilon}$, and if $\epsilon \ne \epsilon^{\prime}$,
then $C^{\epsilon} \ne C^{\epsilon^{\prime}}$.
By Lemma \ref{keyspecial} again, 
for any $\epsilon$ with $\epsilon_1 \le \epsilon \le \epsilon_2$,
there is a sequence of points $\{ d_t^{\epsilon} \} \subset h^{-1}(\beta_j)$
such that $\lim_{t \to \infty} {d_t^{\epsilon} \over \| d_t^{\epsilon} \| }
= C^{\epsilon}$.
This implies that $\dim (D(A_1) \cap D(A_2)) \ge 1$,
which contradicts our assumption.
Thus the limit points are the same point.
\end{proof}

Note that $R_j \in \{ P_1, \cdots , P_a \}$ for $1 \le j \le a + 1$.
Therefore there are $u$, $v$ with $1 \le u, \ v \le a + 1$
and $u \ne v$ such that $R_u = R_v$.
On the other hand, there is $C_1 > 0$ such that
$ST_1(\alpha_u;C_1) \cap ST_1(\alpha_v;C_1) = \{ 0 \}$.
By Lemma \ref{sandwichlemma}, there is $C_2 > 0$ such that
$$
ST_1(h^{-1}(\alpha_u);C_2) \cap ST_1(h^{-1}(\alpha_v);C_2) = \{ 0 \}.
$$
This contradicts the fact that $R_u = R_v$.
Thus $\dim (D(B_1) \cap D(B_2)) = 0$.
\end{proof}

It follows from Claims 1, 2, 3  that if $n \le 3$, then we have
$$
\dim (D(A_1) \cap D(A_2)) = \dim (D(B_1) \cap D(B_2)),
$$
namely the directional dimension is preserved by a bi-Lipschitz homeomorphism.
This is enough to give a comprehensive interpretation
for Oka's family.

\medskip

\end{document}